\documentclass{amsart}
\usepackage{amsfonts,amscd}

\newtheorem{claim}{}[section]
\newtheorem{theorem}[claim]{Theorem}
\newtheorem{lemma}[claim]{Lemma}
\newtheorem{proposition}[claim]{Proposition}
\newtheorem{corollary}[claim]{Corollary}

\def\proclaim #1. #2\par{\medbreak
\noindent{\bf#1.\enspace}{\sl#2}\par\medbreak}
\makeatother

\DeclareMathOperator{\Cdb}{{\mathbb C}}
\DeclareMathOperator{\Rdb}{{\mathbb R}}
\DeclareMathOperator{\Ddb}{{\mathbb D}}
\DeclareMathOperator{\Tdb}{{\mathbb T}}

\DeclareMathOperator{\Real}{Re}

\begin{document}

\title[Dual operator systems]{Dual operator systems}

\date{July 25, 2008.}

\author{David P. Blecher}
\address{Department of Mathematics, University of Houston, Houston,
TX
77204-3008}
\email[David P. Blecher]{dblecher@math.uh.edu}
 \author{Bojan Magajna}
\address{Department of Mathematics,
University of Ljubljana, Jadranska 21, Ljubljana 1000, Slovenia}
\email[Bojan Magajna]{Bojan.Magajna@fmf.uni-lj.si}
\thanks{*Blecher is partially supported by a grant from
the National Science Foundation.
Magajna is  partially supported by the Ministry of Science and
Education of Slovenia. }
\thanks{Some results herein were presented at GPOTS 2008 at University
of Cinncinnatti.}

\begin{abstract}  We characterize weak* closed unital vector spaces
of operators on a Hilbert space $H$.  More precisely, we
%DD
first show that an operator system, which is the dual of an operator
space, can be represented completely isometrically and weak*
homeomorphically as a weak* closed operator subsystem
%DD
of $B(H)$. An analogous result is proved for  unital operator
spaces.
%DD
Finally, we give some somewhat surprising examples of dual unital
operator spaces.
\end{abstract}

\maketitle

\let\text=\mbox

\section{Introduction}
%D
The notion of {\em complete positivity} plays a profound role in the
theory of operator algebras, and also in mathematical physics. The
natural setting for this notion is the category of {\em operator
systems}, namely selfadjoint vector spaces of operators on a Hilbert
space $H$ containing the identity operator $I_H$ (see e.g.\
\cite{CE}, \cite[Section 1.3]{BLM}, \cite{ER}, \cite{P}).  The main topic of our
paper is how operator systems behave with respect to duality.  Also,
we will investigate the duality theory of  {\em unital operator spaces} (defined similarly
to operator systems, but dropping the requirement of selfadjointness).
  Unital operator spaces constitute
an important area of investigation too, for example because this
class includes most interesting operator algebras, selfadjoint or
otherwise, and also includes most interesting function spaces.

Let $X$ be an operator system which, as an operator space, is the
dual of some operator space.
%DD
 By a basic result in operator space theory,
$X$ can be represented completely isometrically and weak*
homeomorphically as a weak* closed subspace, say $X^{\prime}$, of
some $B(H)$, but this theorem fails to guarantee that $X^{\prime}$
is also an operator system, or even unital. In this note we will
show the existence of completely isometric and weak* continuous
%DD
representations of $X$ as a weak* closed operator subsystem of
$B(H)$.  Using notation that will be defined more fully later, we
have:

\begin{theorem}\label{th1} If $X$ is an operator system such that $X$ is a dual as an operator
space, then there exists a completely isometric weak* homeomorphism $\Phi$ from $X$ onto
a weak* closed operator subsystem in $B(\mathcal{H})$ for a Hilbert space $\mathcal{H}$,
which maps the identity into the identity and
hence is also a complete order isomorphism.
\end{theorem}

This is the `operator system' variant of
Sakai's fundamental characterization of von Neumann algebras as
precisely the $C^*$-algebras with a predual.
To prove the theorem, we will first show
%D
that the positive cone of the predual
of $X$ is weak* dense in the positive cone of the dual of $X$.
Similarly, the weak* continuous  states on $X$ are  weak* dense in the states
on $X$ (defined below).  This will provide  enough weak* continuous states
on each $M_n(X)$ to
%D
enable us to construct the desired representation by an `upgrade' of
the proof of Ruan's theorem characterizing operator spaces \cite{ER}.

We will also give a simple metric characterization  of dissipative contractions in unital
operator spaces.  This will enable us to prove a variant of Theorem \ref{th1} for unital
operator spaces:
%D
a unital operator space, which is also a dual operator space, can be
represented completely isometrically and weak* homeomorphically as a
weak* closed unital subspace of $B(\mathcal{H})$ for a Hilbert space
$\mathcal{H}$.   Again, the key tool for this is the weak* density
of
%DD
the weak* continuous  states in the set of all states.

In the last section we will consider various examples, some of them
surprising. For example, the extended and the normal Haagerup tensor
products (see e.g.\ \cite[Section 1.6]{BLM}) of unital dual operator
algebras and spaces, are shown to have
% and construct in this context
natural
weak* continuous unital completely isometric representations on a Hilbert space.
As a consequence of this, it will follow for example
 that many of the algebras met in modern `abstract
harmonic analysis', such as the Fourier algebra $B(G)$,
 are weak* closed unital subspaces of a von Neumann algebra,
or unital subspaces of a $C^*$-algebra.

Turning to notation, we write $X^{\sharp}$ for the dual of $X$, and $X_{\sharp}$
for a predual of $X$.  Weak* continuous maps will often be referred to as `normal',
as is customary.  If $X$ has a `unit' $1$, we write $S(X)$ for the set of all states
on $X$, that is, the set of all functionals $\rho\in X^{\sharp}$ with $\rho(1)= \Vert \rho \Vert = 1$.
We write $X^h$ for the hermitian elements of $X$, namely the elements $x$ with
$\rho(x) \in \Rdb$ for all states $\rho$
%D
 (see \cite{BD,V}).
 %DD  Bojan, is this what you wanted?  ADD
 If $X$ is an operator system then we write $X^+$ for the positive cone of $X$.
We also write $X^{{\sharp}+}$ (resp.\ $X_{\sharp}^+$) for the
positive cone of $X^{\sharp}$ (resp.\ $X_{\sharp}$): these are the
functionals (resp.\ normal functionals) satisfying $\rho(x)\geq0\
\forall x\in X^+$. In this case the state space $S(X)$ consists of
all $\rho\in X^{{\sharp}+}$ with $\rho(1)=1$. The hermitian elements
also have several interesting alternative characterizations in this
case, for example $x \in X^h$ iff $\|1+itx\|^2=1+t^2\|x\|^2$ for all
$t\in\mathbb{R}$, and iff $\|1+itx\| = 1 + o(t)$ (see e.g.\
\cite{BD,V}).

%DD
An abstract characterization of unital operator spaces may be found
in \cite{BN}.  The reader may also find metric characterizations of
operator systems there,
%DD
and the fact that the involution on a dual operator system is weak*
continuous. The famous order theoretic characterization of operator
systems may be found in \cite{CE,P}. By a {\em subsystem} of an
operator system $X$ we mean a selfadjoint subspace containing the
identity of $X$. A {\em complete order isomorphism} is a linear
isomorphism $T$ between operator systems such that $T$ and $T^{-1}$
are completely positive (and thus necessarily preserve the
involution).
 If $T$ is unital (that is, $T(1) = 1$), then
this turns out to be equivalent to $T$ being a surjective complete isometry.
 For this and any additional details or definitions concerning
 operator spaces and systems the reader is referred to the texts
cited in our bibliography.

\section{Weak* density of normal states and dissipative elements}

In the proofs below we will use,  without explicitly mentioning it,
the theorem asserting that a convex subset of a dual Banach space $V$ is weak* closed if and only if its intersection with each closed
ball in $V$ is weak* closed.
This is known as the Krein--Smulian theorem
in operator  theory, but in \cite[4.44]{F} it is attributed to Banach and  Dieudonn\'{e}.

The fact that $X_{\sharp}^+$ is weak* dense in $X^{{\sharp}+}$ turns out to be a consequence of the fact that
the cone $X^+$ of positive elements of $X$ is weak* closed.

\begin{lemma}\label{le1}
%DD
Let $X$ be an operator system which is a dual Banach space.
Then \begin{itemize}
\item [(i)] $X^+$ is weak* closed.

\item [(ii)] $X^+$ is the dual cone to $X_{\sharp}^+$.  That is, denoting $(X_{\sharp}^+)^{\circ}=\{x\in X:
\langle x,\omega\rangle\geq0\ \forall\omega\in X_{\sharp}^+\}$, we have that $(X_{\sharp}^+)^{\circ}=X^+$.

\item [(iii)] $X_{\sharp}^+$ is weak* dense in $X^{{\sharp}+}$.

\item [(iv)] $X_{\sharp}\cap S(X)$ is weak* dense in $S(X)$.
\end{itemize}
\end{lemma}

\begin{proof} (i) \ By  \cite[3.7]{BN}, the set $X^h$ of all hermitian elements in $X$ is
weak* closed. Denoting by $B_X$ the unit ball of $X$, we have that an element $x\in X^h\cap
B_X$ is in $X^+$ if and only if $\|1-x\|\leq1$. From this we deduce that $X^+\cap B_X$ is
weak* closed, hence  $X^+$ is weak* closed (since
$X^+$ is convex and $X^+\cap rB_X=r(X^+\cap B_X)$ for all $r>0$).

(ii)  \ Suppose that there exists a selfadjoint $x\in
(X_{\sharp}^+)^{\circ}\setminus X^+$. Then (since $X^+$ is weak*
closed and convex) by a geometric form of the Hahn--Banach theorem
there exist an $\alpha\in\mathbb{R}$ and an $\omega\in X_{\sharp}$
such that $\Real \, \langle y,\omega\rangle\leq\alpha$ for all $y\in
X^+$ and $\Real \, \langle x,\omega\rangle >\alpha$. We may replace
$\omega$ by $\frac{1}{2}(\omega+\omega^*)$, where
$\omega^*(v):=\overline{\omega(v^*)}$ for all $v\in X$. This is a
weak* continuous functional, since the involution is weak*
continuous on $X$ by \cite{BN}. Since $X^+$ is a cone, we may take
$\alpha=0$. Then $\langle y,\omega\rangle\leq0$ for all $y\in X^+$,
while $\langle x,\omega\rangle >0$. But the first relation means
that $-\omega\in X_{\sharp}^+$, hence $\langle
x,-\omega\rangle\geq0$ since $x\in(X_{\sharp}^+)^{\circ}$, which
contradicts the second relation.  This is all that is needed for the
proof of (iii) below.  For general $x \in (X_{\sharp}^+)^{\circ}$,
by (iii) and  weak* density we have $\varphi(x) \geq 0$ for all
$\varphi \in X^{{\sharp}+}$, so $x \in X^+$.

(iii) \ Suppose that there exists a $\rho\in
X^{{\sharp}+}\setminus\overline{X_{\sharp}^+}^{w*}$. By a geometric
form of the Hahn--Banach theorem, similarly to the proof of (ii),
there exists an $x\in X^h$ such that  $\langle x,\omega\rangle\leq
0$ for all $\omega\in X_{\sharp}^+$ and $\langle x,\rho\rangle >0$.
The first relation implies that $-x\in(X_{\sharp}^+)^{\circ}$, and
hence $-x\in X^+$ by (ii).  Consequently $\langle
-x,\rho\rangle\geq0$ since $\rho\in X^{{\sharp}+}$, which
contradicts the second relation.

(iv) \ Let $\rho\in S(X)$. Since $\rho\in X^{{\sharp}+}$, by (iii) there exists a net
$\psi_{\alpha}\in X_{\sharp}^+$  weak* converging to $\rho$. In particular $\psi_{\alpha}(1)\to
\rho(1)=1$, hence the net of states $\phi_{\alpha}:=\psi_{\alpha}(1)^{-1}\psi_{\alpha}\in X_{\sharp}
\cap S(X)$ also weak* converges to $\rho$.
\end{proof}

{\bf Remark.} \
%DD
Parts (iii) and (iv) of the above lemma can be
proved via a result from the theory of order unit and base normed
spaces (see e.g.\ \cite{AE}), but we preferred the more `elementary'
proof given
 above.
Indeed since the selfadjoint part of $X$ is an order unit space, it
follows
%DD
that the selfadjoint part of $X_{\sharp}$ is a base normed
space, and hence it has all the useful properties of
%DD
such spaces.
%base normed spaces.  For example, it follows that every
%selfadjoint contractive normal functional on $X$ is a norm limit of
%functionals of the form $t \varphi - (1-t) \psi$ for normal states
%$\varphi, \psi$ on $X$ and $t \in [0,1]$.

\smallskip

\begin{proof}[Proof of Theorem {\rm  \ref{th1}}] Let $n\in\mathbb{N}$, set $Y=M_n(X)$, let $x^*=x\in Y$
be a fixed element with $\|x\|=1$, and choose a state $\rho$ on $Y$ such that $|\rho(x)|=
\|x\|=1$. Given $\varepsilon>0$, by Lemma \ref{le1} (iv) there exists a state
$\phi\in Y_{\sharp}\cap S(Y)$ such that $|\phi(x)|>1-\varepsilon$.

Put $A=M_n(\mathbb{C})$. Inspired by the proof of Ruan's theorem \cite[p. 30--34]{ER}, we
{\em claim} that there exists a state $\omega$ on $A$
such that for all $m\in\mathbb{N}$
\begin{equation}\label{2}|\phi(avb)|\leq\omega(aa^*)^{1/2}\omega(b^*b)^{1/2}\ , \qquad
v\in M_{nm}(X)\ , \|v\|\leq1,\ a,b^*\in M_{1,m}(A).
\end{equation}

Assuming the claim (which will be proved below), denoting by $\pi$
the cyclic representation on the Hilbert space $H$ constructed from
$\omega$ (necessarily equal to the multiple $k\cdot{\rm id}$ of the
identity representation up to a unitary equivalence, $k\leq n$,
$H\subseteq(\mathbb{C}^n)^n$) and by $\xi\in H$ the corresponding
unit cyclic vector, (\ref{2}) may be rewritten:
$$|\phi(avb)|\leq\|\pi_{1,m}(a)^*\xi\|\|\pi_{m,1}(b)\xi\|\ , \qquad  a,b^*\in M_{1,m}(A), v \in M_{nm}(X)\ . \|v\|\leq1, $$
Taking $m=1$ we see that for each fixed contraction $v\in M_{n}(X)$,
the map $$(\pi(b)\xi, \pi(a^*)\xi)\mapsto\phi(avb)$$ is a
contractive sesquilinear form on $[\pi(A)\xi]
\times[\pi(A)\xi]=H\times H$.  Hence there is a contraction $S(v)$
on $H$ such that
\begin{equation}\label{3}\phi(avb)=\langle S(v)\pi(b)\xi,\pi(a^*)\xi\rangle\
, \qquad \ a,b\in A.
\end{equation}
We may extend $S$ to a linear map $S:M_{n}(X)\to B(H)$ satisfying (\ref{3}) for all $v\in M_n(X)$. Since $\xi$ is
cyclic for $\pi(A)$, it follows from (\ref{3}) that $S$ is a bimodule map over
$A=M_n(\mathbb{C})$.  That is
$$S(cvd)=\pi(c)S(v)\pi(d)\ , \qquad v\in M_{n}(X),\ c,d\in A.$$
Thus, $S$ must be the amplification
$T_n$ of a map $T:X\to B(\mathbb{C}^k)$, and $S$ (hence $T$) is automatically
completely contractive by \cite{S} since $\pi$ is cyclic.

Further,  $S$ (hence $T$) is unital. To show this, first note that
$1=\phi(1)=\langle S(1)\xi,\xi\rangle$ (by (\ref{3})), which implies
(since $\|S(1)\|\leq1$ and $\|\xi\|=1$) that $S(1)\xi=\xi$. The
bimodule property implies that $\pi(c)S(1)= S(c1)=S(1c)=S(1)\pi(c)$
for all $c\in A$. Hence $S(1)\pi(c)\xi=\pi(c)S(1)\xi=\pi(c)\xi$,
implying that $S(1)=1$ since $\xi$ is cyclic for $\pi(A)$.

Moreover, (\ref{3}) implies that
\begin{equation}\label{1}\phi(avb)=\langle\pi(a)T_n(v)\pi(b)\xi,\xi\rangle\
, \qquad  v\in M_n(X),\ a,b\in A.
\end{equation}
By a standard argument, it follows that $T_n$ (hence $T$) is weak*
continuous (since $\phi$ is weak* continuous and $\xi$ is cyclic for
$\pi(A)$,  and using also the theorem mentioned in the first
paragraph of our paper). Since $|\langle T_n(x)\xi,\xi\rangle|=
|\phi(x)|>1-\varepsilon$ implies that $\|T_n(x)\|>1-\varepsilon$,
taking the direct sum of all such maps $T$ (over all selfadjoint $x$
in the unit sphere of $M_n(X)$, all $n\in\mathbb{N}$ and
$\varepsilon\in (0,1/2)$), we obtain a weak* continuous unital
complete isometry $\Phi$ from $X$ into $B(\mathcal{H})$ for a
Hilbert space $\mathcal{H}$. Then $\Phi(X)$ is automatically weak*
closed and $\Phi$ is a weak* homeomorphism onto its range (this is a
well known consequence of the weak* compactness of closed balls and
the result mentioned in the beginning of the section). By facts
mentioned in the introduction,
%D, for example, the metric characterization of hermitian elements mentioned at the end of Section 1,
%($y$ is hermitian iff $\|1+ity\|^2=1+t^2\|y\|^2$ for
%%all $t\in\mathbb{R}$ \cite{V}),
%a unital complete isometry $\Phi$ necessarily maps hermitian elements to
%hermitian, hence preserves the
$\Phi$ is a complete order isomorphism, and preserves the involution.

It remains to prove the claim, that is, to show the existence of a state $\omega$ on $A$
satisfying (\ref{2}).  For this, it suffices to show that there exists a state $\omega$ on $A$
such that  for every $m\in\mathbb{N}$ we have
\begin{equation}\label{4}\phi(cvc^*)\leq\omega(cc^*),\ \mbox{where}\ c\in M_{1,m}(A),\
 v\in M_{nm}(X),\  v=v^*,\ \|v\|\leq1.
\end{equation}
Namely, using (\ref{4}) with $c=a+b^*$ and with $c=a-b^*$, and the fact that
$\phi^*=\phi$ and $v^*=v$, we obtain
$$\begin{array}{lll}\Real \, \phi(avb)&=&\frac{1}{2}\phi(avb+b^*va^*)\\
 &=&\frac{1}{4}\phi((a+b^*)v(a+b^*)^*-(a-b^*)v(a-b^*)^*)\\
 &\leq&\frac{1}{4}[\omega((a+b^*)(a+b^*)^*)+\omega((a-b^*)(a-b^*)^*)]\\
 &=&\frac{1}{2}\omega(aa^*+b^*b).
 \end{array}$$
 Replacing  $a$ by $ta$ and $b$ by $t^{-1}b$
in this estimate, and taking the infimum over
 all $t\in(0,\infty)$, we get
 \begin{equation}\label{5}\Real \, \phi(avb)\leq\omega(aa^*)^{1/2}\omega(b^*b)^{1/2}
 \end{equation}
 for all selfadjoint contractions $v$. However, for a general $v$ we may write
$$avb=[a,0]\left[\begin{array}{ll}
0&v\\
v^*&0
\end{array}\right]\left[\begin{array}{l}
0\\
b
\end{array}\right],$$
hence (\ref{5}) must hold in general. Replacing $v$ by a suitable multiple
$zv$, where $|z|=1$, (\ref{2}) follows.

Finally, to prove the existence of a state $\omega$ satisfying (\ref{4}), consider the cone $C$ consisting of all functions $f_{c,v}$
on the state space of $A$ of the form $f_{c,v}(\omega)=\omega(cc^*)-\phi(cvc^*)$, where
$c\in M_{1,m}(A)$ (for all $m\in\mathbb{N}$) and $v\in M_{nm}(X)$ is selfadjoint with
$\|v\|\leq1$. A standard argument shows that $C$ is indeed a cone and that each function in
$C$ takes a nonnegative value at some point of the state space $S(A)$ of $A$.
Thus by \cite[2.3.1]{ER}
there exists a point $\omega\in S(A)$ at which all the functions of $C$ are nonnegative, which
proves (\ref{4}).
\end{proof}

Using Lemma \ref{le1} we can also prove a version of Theorem \ref{th1} for function systems.
\begin{corollary}
 \label{dfsy}  Let $X$ be a function system
(by which we mean a subsystem of $C(K)$ for $K$
compact) with a Banach space
 predual.  Then there exists a weak*-homeomorphic unital order isomorphism
 and isometry of $X$ onto a weak* closed subsystem of a commutative
 $W^*$-algebra.  The  $W^*$-algebra may be taken to be
 $\ell^\infty(\Omega)$ where $\Omega$ is the set of normal states of
 $X$.  \end{corollary}

\begin{proof}   Let $\Omega$ be as above, and let $\theta : X \to
\ell^\infty(\Omega)$ be the canonical contraction.  This is easily
seen to be weak* continuous.   Given $x \in X$, let $\varphi$ be a
state of $X$ with $|\varphi(x)| = \Vert x \Vert$ (this is possible
by the definition of a function system).  Let $(\varphi_t)$ be a net
of normal  states on $X$ converging weak* to $\varphi$. Since
$|\varphi_t(x)| \to |\varphi(x)|$, it follows that $\theta$ is an
isometry.  The proof is now completed by the general functional
analytic principles in the paragraph before   Equation (\ref{4}).
\end{proof}

%{\bf Remark.}
Comparing Theorem \ref{th1} to Sakai's characterization of
von Neumann algebras mentioned in the introduction, it is natural to ask
if there exists an operator system $X$ with a {\em Banach space predual},
%for which the involution is even weak* continuous,
but for which there exists no   weak*-homeomorphic unital complete
isometry of $X$ onto a weak* closed subsystem of a $W^*$-algebra?
Indeed, an example of such was exhibited in \cite[Proposition 2.1]{BM}.
There the involution was even weak* continuous, but $X$ had no
operator space predual.

The following result settles another natural question
that arises when considering the famous `unicity of predual' of von Neumann algebras:

\begin{proposition}   An operator system may have more than
%DD one operator space
predual.
\end{proposition}

\begin{proof}   This follows by a routine technique (see e.g.\ \cite[Corollary
2.7.8]{BLM}):
 Let $X$ be an operator space
with two distinct operator space preduals.  With respect to each of
these two preduals, there are two complete isometries from $X$ onto
weak* closed subspaces of $B(H)$ and of $B(K)$ respectively.  The
corresponding `Paulsen systems' in $M_2(B(H))$ and $M_2(B(K))$ are
weak* closed operator systems, and are completely isometric to each
other via a unital complete order isomorphism by `Paulsen's lemma'
(e.g.\ in the form in \cite[Lemma 1.3.15]{BLM}). However if this
isomorphism was a weak*-homeomorphism, then so is its restriction to
the copies of $X$, which gives a contradiction as in e.g.\
\cite[Corollary 2.7.8]{BLM}.
\end{proof}

%\begin{proposition}   {\bf Remark.}
  It is easy to see that the normal states span the predual.  Indeed:

\begin{proposition}  If
%DD
$X$ is a weak* closed unital operator space then the
any normal functional on $X$ of norm $< 1$ may be written as
$\varphi_1 + i \varphi_2$, where each of $\varphi_1, \varphi_2$ is
of the form
%DD
$s \xi - t \psi$ for normal states $\xi, \psi$ on $X$ and $s, t \in
[0,1], s + t < 1$.
\end{proposition}

\begin{proof}
%DD
If $X$ is weak* closed in $B(H)$, then
 $X_{\sharp}$ is a quotient of
$B(H)_{\sharp}$.   Thus we may assume that $X = B(H)$, by
considering preimages of functionals.   However the von Neumann
algebra case of the result is well known.
\end{proof}

To prove a version of Theorem \ref{th1} for unital operator spaces,
 we will need a suitable variant of Lemma \ref{le1}. Since there are
 usually not enough hermitian elements in such spaces, we will consider dissipative elements instead.

\begin{lemma}\label{le2} An element $x$ with $\|x\|\leq1$ in a unital operator space $X$ is dissipative (has numerical range contained
in the closed left half-plane) if and only if
\begin{equation}\label{21}\|1+tx\|^2\leq1+t^2\ \mbox{for all}\ t\in[0,1].\end{equation}
Consequently, the set $D_X$ of all dissipative elements in a unital  operator space $X$ with
a Banach space predual is weak*  closed.
\end{lemma}

\begin{proof} Suppose that $\|x\|\leq1$.  If $x$ is dissipative, then for all $t\in[0,1]$
$$\|1+tx\|^2=\|1+2t\Real \, x+t^2x^*x\|\leq\|1+2t\Real \, x\|+t^2\leq1+t^2.$$
Conversely, if the numerical range $W(x)$ of $x$ contains some $\alpha$ with $s:=\Real \, \alpha>0$,
then $s\in W(\Real \, x)$, and hence
$$\|1+tx\|^2\geq\|1+2t\Real \, x\|-t^2\|x\|^2\geq1+2ts-t^2>1+t^2,\ \mbox{if}\ t\in(0,s).$$
It follows that the intersection of $D_X$ with the unit ball of $X$
is weak* closed, hence (since $D_X$ is convex and closed under
multiplication by positive scalars) the same holds for any ball
around $0$.   Consequently $D_X$ is weak* closed.
\end{proof}

\begin{corollary}\label{co1} If $X$ is a unital operator space with a Banach space predual,
then the set
%DD $X_{\sharp}\cap S(X)$
of all
weak* continuous states on $X$ is weak* dense in S(X).
\end{corollary}

\begin{proof} Denote by $X^{{\sharp}+}$ the set of all nonnegative multiples of states on $X$
and by $(D_X)^{\circ}$ the set of all functionals $\rho\in X^{\sharp}$ such that
$\Real \, \rho(x)\leq0$ for all $x\in D_X$. Clearly $X^{\sharp+}\subseteq (D_X)^{\circ}$, but we claim
that the two sets are equal. To prove this, let $\rho\in (D_X)^{\circ}$.
Note that $it1\in D_X$ for all $t\in\Rdb$ implies that  $t\Real \, (i\rho(1))\leq0$, hence
$\rho(1)\in\Rdb$.  Since $-1\in D_X$, $\rho(1)\geq0$.
Since $x-t1\in D_X$ for each $x\in X$ and $t\geq\|x\|$, we have that $\Real \, \rho(x)\leq
\|x\|\rho(1)$. Replacing in this inequality $x$ by $zx$ for all $z\in\Cdb$ with $|z|=1$, it follows that
$|\rho(x)|\leq\|x\|\rho(1)$, hence $\rho\in X^{\sharp+}$.
Set $X_{\sharp}^+=X_{\sharp}\cap X^{{\sharp}+}$ and $(D_X)_{\circ}=(D_X)^{\circ}\cap X_{\sharp}$.
Then $X_{\sharp}^+=(D_X)_{\circ}$.

From now on the proof is similar to the proof of Lemma \ref{le1}.
Since $D_X$ is weak* closed by Lemma \ref{le2}, a bipolar type
argument (as in the proof of Lemma \ref{le1}) shows that
$D_X=((D_X)_{\circ})^{\circ}$ (where the polar $A^{\circ}$ of a
subset $A$ of $X_{\sharp}$ is defined as $A^{\circ}=\{x\in X:\Real
\, \langle a,x\rangle\leq0\ \forall a\in A\}$) and that
$(D_X)_{\circ}$ is weak* dense in $(D_X)^{\circ}$. By the previous
paragraph this means that $X_{\sharp}^+$ is weak* dense in
$X^{\sharp+}$ and it follows (as in the proof of Lemma
\ref{le1}(iv)) that $X_{\sharp}\cap S(X)$ is weak* dense in $S(X)$.
\end{proof}

Similarly to Corollary \ref{dfsy}, we now have:

\begin{corollary} A unital function space $X$ with a Banach space
predual is weak*-homeomorphic, order isomorphic, and isometric (via
the same unital map), with a weak* closed unital subspace of
$\ell^\infty(\Omega)$, where $\Omega$ is the set of normal states of
$X$.  \end{corollary}

%D
We will study a dual unital operator space $X$ in terms of its
`canonical operator system' $X + X^*$.  That is, we may assume that
$X\subseteq B(H)$ for a Hilbert space $H$ as a unital operator
space, and consider the subsystem $X+X^* = \{
%DD
x + y^* : x,y \in X \}$ in $B(H)$.  We use some of Arveson's results
concerning $X + X^*$ (see \cite{Arv} or 1.3.6 and 1.3.7 in
\cite{BLM}), such as that this system is independent of $H$. If $X$
is also a dual operator space, then one might at first hope that
%there exists a representation of $X$ for which
$X + X^*$ is also a dual space containing $X$ as a weak* closed
subspace, and we could directly apply Theorem \ref{th1}. However
this is false.  Indeed we recall the result from Banach space theory
\cite{F}
 that for weak* closed subspaces $E$ and $F$ of a dual space,
$E + F$ is weak* closed if and only if it is norm closed (i.e.\
complete). Thus $X + X^*$ is weak* closed with respect to one (or
every) unital weak* continuous completely isometric representation
of $X$
%DD
if and only if $X + X^*$ is norm complete.
%DD
If, for example,  $H^\infty(\Ddb) + H^\infty(\Ddb)^*$ were a dual
space, then it is weak* closed in $L^\infty(\Tdb)$. Since
$H^\infty(\Ddb) + H^\infty(\Ddb)^*$ is weak* dense in
$L^\infty(\Tdb)$, we would have $H^\infty(\Ddb) + H^\infty(\Ddb)^* =
L^\infty(\Tdb)$, which is false (see the discussion on p.\ 181--182
of \cite{H}).
%Thus we will have to work harder.

\begin{theorem}\label{th2}
%DD
If $X$ is  unital operator space, which is also a dual operator
space, then there is a completely isometric unital weak*
homeomorphism from $X$ onto  a weak* closed unital subspace of
$B(\mathcal{H})$ for a Hilbert space $\mathcal{H}$.
\end{theorem}

\begin{proof}  Set $Y = X + X^*$ as discussed above.
%We may assume that $X\subseteq B(H)$ (for a Hilbert space $H$) as a unital operator space and denote by
Extend each $\phi\in X_{\sharp}\cap S(X)$ to a state $\tilde{\phi}$
on $Y$. Since $\tilde{\phi}(x+y^*)=\phi(x)+\overline{\phi(y)}$ for
all $x,y\in X$, it follows from Corollary \ref{co1} that the set
$\tilde{S}$ of all such extensions is weak* dense in $S(Y)$, and
similarly on all matrix levels $M_n(Y)=M_n(X)+M_n(X^*)$. Now, using
the states in $\tilde{S}$ for constructing representations of $Y$ as
in the proof of Theorem \ref{th1}, we get a
%DD
unital complete isometry
$\Phi$ from $Y$ onto a subspace of some $B(\mathcal{H})$ such that
$\Phi|X$ is weak* continuous, hence a weak* homeomorphism from $X$
onto a weak* closed subspace of $B(\mathcal{H})$.
\end{proof}

{\bf Remark.}
Every unital dual operator space $X$ has a canonical
`dual operator system envelope'
$\tilde{X}$.  That is, $X$ is weak* closed in $\tilde{X}$,
 $\tilde{X}$ is the weak* closure of $X + X^*$, and for any weak* continuous
 complete isometry
$u : X \to B(H)$ there exists a weak* continuous completely positive
unital
%DD
map $\Phi$ from $\tilde{X}$ onto the weak* closure of $X + X^*$ in
$B(H)$,
%DD
which extends $T$.  This may be constructed in a standard
way by setting ${\mathcal F}$ be the set of all (or enough) weak*
continuous
%DD
unital complete isometries $T : X \to B(H)$, where $H$ may be any
Hilbert space (of cardinality small enough).  Write $H$ here as
$H_T$, and let $i(x) = (T(x))_{T \in {\mathcal F}} \in B(\oplus_{T
\in {\mathcal F}} \; H_T)$, for $x \in X$. Set $\tilde{X}$ equal to
the weak* closure of $i(X) + i(X)^*$, this will clearly have the
desired property.
%DD
It would be interesting to investigate the spaces $X$ for which
$\Phi$ above is always completely isometric.  Equivalently, for
which $X$ is the weak* closure of $T(X) + T(X)^*$, as an operator
system, independent of the representation $T$ of $X$.
%DD
If this were always true for dual unital function spaces for
example, then it is easy to see that every normal state on a weak*
closed unital subspace $X$ of $L^\infty(\Omega)$ extends to a normal
state on the weak* closure of $X + X^*$ in $L^\infty(\Omega)$.  This
seems to be known to usually be false (see e.g.\ \cite{HR}).
%is related to known variants of the
%Gleason-Whitney theorem in function theory, a well understood area
%of .
%DD Along these lines, we remark that
%the argument just above Theorem \ref{th2} shows that if
%$X + X^*$ is weak* closed in one weak* continuous completely isometric representation of
%$X$, then it is weak* closed in every other such representation.

%D
\section{Examples}
%Representation of the extended and the normal Haagerup tensor product of unital operator spaces}

First, we indicate how the normal Haagerup tensor product
  of dual
unital operator algebras can be represented faithfully and weak*
homeomorphically into $B(H)$.

If $M$ is a von Neumann algebra, then we recall from \cite{EK} that
$M\otimes ^{\sigma {\rm h}} M \cong CB_{M'}(B(H))$ isometrically,
and weak* homeomorphically. In fact it is not difficult to show that
this is a complete isometry (see \cite[4.4]{DANP} for a more general
result).
%DD
As usual, $CB$ denotes the completely bounded maps.

If  $M$ and $N$ are von Neumann algebras, then we will write $M * N$
for the `universal von Neumann algebra amalgamated free product',
amalgamating over the scalars. This is a von Neumann algebra with
the appropriate  free product universal property, appropriate to
`extending to $M * N$ every pair of normal $*$-representations of
$M$ and $N$ on a common Hilbert space'.  This universal property
forces the object to be unique. That it exists may be proved by a
routine `soft' argument using large direct sums (somewhat in the
spirit of  \cite[3.1]{B} or  the last Remark), or by taking the
direct sum of all cyclic representations $\pi$ of the
$C^*$-algebraic free product of $M$ and $N$ (see \cite[p. 98]{Pi})
such that $\pi|M$ and $\pi|N$ are normal.

%, or just take the von Neumann envelope of the amalgameted free product
%of the C$^*$-algebras $M$ and $N$ \cite[p. 98]{Pi}.

\begin{theorem}  \label{win}
\begin{itemize} \item [(1)] If  $M$ and $N$ are unital dual operator  algebras, then $M\otimes^{\sigma {\rm h}} N$
may be identified (completely isometrically and weak* homeomorphically), with
the weak* closure of $M \otimes N$ in $W^*_{\rm max}(M) \otimes^{\sigma {\rm h}} W^*_{\rm max}(N)$.
Here  $W^*_{\rm max}(\cdot)$ is the `maximal  von Neumann algebra' of \cite{BS}.

\item [(2)] If  $M$ and $N$ are von Neumann algebras, then
 $M\otimes^{\sigma {\rm h}} N$
may be identified (completely isometrically and weak* homeomorphically), with
the weak* closure of $M  N$ in the universal  von Neumann algebra amalgamated free
product $M * N$.

\item [(3)] If  $M$ and $N$ are unital dual operator  algebras, then $M\otimes^{\sigma {\rm h}} N$
is a Banach algebra with product $(a \otimes b)(a' \otimes b') = aa'  \otimes b' b$.

\item [(4)] If  $M$ and $N$ are unital dual operator  algebras, then $M\otimes^{\sigma {\rm h}} N$
is (completely isometrically and weak* homeomorphically) a weak* closed
unital subspace of some von Neumann algebra.
\end{itemize}
\end{theorem}

\begin{proof}  (1) \ To prove this, let $u : M \times N \to B(H)$ be separately weak* continuous and completely
contractive. As Effros and Ruan showed (see e.g.\
\cite[1.6.10]{BLM}), we may write $u(x,y) = R(x) S(y)$ for weak*
continuous   complete contractions $R, S$.  By \cite[Theorem
2.7.10]{BLM} we may rewrite this as $u(x,y) = a \pi(x) b \rho(y) c$
for unital weak* continuous   completely contractive homomorphisms
$\pi, \rho$, and contractive operators $a, b, c$.   As in the proof
of e.g.\
%DD
\cite[Theorem 5.13]{Pi}, we may assume that $b = 1$. Let
$\tilde{\pi}, \tilde{\rho}$ be the canonical extensions to $W^*_{\rm
max}(M)$ and $W^*_{\rm max}(N)$ respectively.  The map $a
\tilde{\pi}(x) \tilde{\rho}(y) c$ on $W^*_{\rm max}(M) \times
W^*_{\rm max}(N)$ induces a weak* continuous   complete contraction
on $W^*_{\rm max}(M)  \otimes^{\sigma {\rm h}} W^*_{\rm max}(N)$.
Let $\tilde{u}$ be the restriction of the latter map to the weak*
closure $E$ of $M \otimes N$.  Clearly $\tilde{u}(x \otimes y) =
u(x,y)$ for $x \in M, y \in N$.  Hence $E$ has the universal
property of $M \otimes^{\sigma {\rm h}} N$, and thus it follows that
the canonical map $M \otimes^{\sigma {\rm h}} N \to E$ is a weak*
continuous completely isometric surjection.

(2) \ This is proved analogously to the same fact for the Haagerup
tensor product
%DD
(see \cite[Theorem 5.13]{Pi}).  One only needs the
free product universal property,
 and the methods in the previous paragraph to show that
the appropriate subspace of  $M * N$ has the universal
property of $M \otimes^{\sigma {\rm h}} N$.

(3) \ By (1) we may  suppose that $M, N$ are von Neumann algebras.
Also we may  suppose that $M = N$ by the trick of letting
 $R=M\oplus N$.  It is easy to argue that
$M\otimes ^{\sigma {\rm h}}N \subset R\otimes ^{\sigma {\rm h}}R$, since $M$ and $N$ are
appropriately complemented in $R$.  However, by \cite{EK} we have
$M\otimes ^{\sigma {\rm h}} M \cong CB_{M'}(B(H))$, a Banach algebra.

(4)  \ Follows by combining (1) and (2).
\end{proof}

\begin{corollary} \label{cb}  For any Hilbert space $H$,
the space $CB(B(H))$
%(resp.\   $NCB(B(H))$)
%DD
 is
a weak* closed unital subspace of a von Neumann algebra.
\end{corollary}

\begin{proof}
Apply the above mentioned complete isometry
%DD
from \cite{EK}, and Theorem \ref{win} (4).
\end{proof}

It follows from the last corollary,
that many of the algebras met in modern `abstract
harmonic analysis', such as the Fourier algebra $B(G)$,
 are weak* closed unital subspaces of a von Neumann algebra,
or unital subspaces of a $C^*$-algebra. Indeed there has been quite
a lot of work in recent years, in papers by M. Neufang and
coauthors, and others (see e.g.\ \cite{JNR} and references therein),
showing that many such algebras are completely isometric to unital
subalgebras of
%DD $NCB(B(H))$ or
$CB(B(H))$.
We are indebted to Matthias Neufang for discussions on the state of this topic.

\begin{corollary} \label{cb2}  For any Hilbert space $H$,
the `completely bounded norm' of the row $[ Id  \; \;  T  ]$ in $R_2(CB(B(H))$  is
%DD
 $\sqrt{1 + \Vert T \Vert_{\rm cb}^2}$ for any  $T \in
CB(B(H))$.  Here $Id$ is the identity map on $B(H)$.
%DD
Thus  for any
$T \in CB(B(H))$
 of norm $1$, one can find a matrix $[x_{ij}]$ in the unit ball of
$M_n(B(H))$ with
the norm of   $[ [x_{ij}]  \; \;  [T(x_{ij})]  ]$  arbitrarily
close to $\sqrt{2}$.  Similarly with $R_2$ replaced by
$C_2$.
\end{corollary}

%DD
{\bf Remark.}  Note that $CB(X)$ is not in general a unital operator
space, even if $X$ is finite dimensional.  For example, take a
unital finite dimensional Banach algebra $A$ which is not unitally
isometric to a unital subspace of a $C^*$-algebra (for example there
exist 4 dimensional Banach algebras of numerical index
$\frac{1}{e}$, see e.g.\ \cite[p.\ 112]{BD}, whereas any unital
operator space is easily seen to have index $\geq \frac{1}{2}$). The
assertion follows since the matrix normed algebra Max$(A)$ is a
subspace of $CB({\rm Max}(A))$ via the regular representation.

 On the other hand, if  $M,N\subseteq B(H)$
are von Neumann algebras containing a von Neumann subalgebra $A$ in
their intersection, then the space of all
$M^{\prime},N^{\prime}$-bimodule maps
$CB_{M^{\prime}}(A^{\prime},B(H))_{N^{\prime}}$ is a weak* closed
unital subspace of a von Neumann algebra.  Indeed, this space is
completely isometric and weak* homeomorphic to $M \otimes^{\sigma
h}_A N$ by \cite[4.4]{DANP}. However using ideas of Ozawa, one can
show that $M
 \otimes^{\sigma h}_A N$ may be embedded into the (maximal) von
Neumann amalgamated free product $M \; *_A \; N$.  We will perhaps
present the details elsewhere.

\begin{theorem}  \label{ng}  If $X, Y$ are  unital operator spaces,
then so is $X \otimes_{\rm h} Y$.
%  If $X, Y$ are  unital dual operator algebras,
%then $X \otimes^{\sigma {\rm h}} Y$ and  $X \otimes_{\rm w^*  h} Y$
%are unital operator spaces.
Indeed $X \otimes_{\rm w^*  h} Y$ is a
weak* closed unital-subspace of a von Neumann algebra,
 if $X, Y$ are both weak* closed unital-subspaces of  von Neumann algebras.
%unital dual operator spaces.
\end{theorem}

\begin{proof}   If $X, Y$ are subspaces of $C^*$-algebras $A, B$, then
$X \otimes_{\rm h} Y \subset A \otimes_{\rm h} B$, and it follows
from \cite[Theorem 5.13]{Pi} that $X \otimes_{\rm h} Y$ is a unital-subspace of
the full amalgamated free product $C^*$-algebra $A \, *_{\Cdb} \, B$.

%DD Text supressed:
%IGNORE THIS PARAGR FOR NOW: There is a universal unital dual operator algebra $WOA(X)$
% associated with any unital dual operator space $X$,
%with the following universal property:  given any unital weak* continuous complete
%contraction $T : X \to B(H)$, there is a unital weak* continuous completely
%contractive homomorphism $WOA(X) \to B(H)$ extending $T$.  This may be proved
%by the usual soft analysis trick.    We deduce that
%$X \otimes^{\sigma {\rm h}} Y \subset WOA(X) \otimes^{\sigma {\rm h}} WOA(Y)$.  Indeed
%given any separately weak* continuous $u : X \times Y \to B(H)$,
%we may write $u(x,y) = R(x) S(y)$ for   Now $R : X \to B(K,H)$
%is a weak* continuous complete contraction.
%If $X, Y$ are  unital operator spaces,
%then  $X \otimes_{\rm h} Y$ is a unital-subspace of
%$C^*(X) \otimes^{\rm h} C^*(Y)$, which in turn is a unital-subspace of the full
%$C^*$-algebra amalgamated free product of $C^*(X)$ and $C^*(Y)$.

%From Theorem \ref{win} (4),  if $M, N$ are  unital dual operator algebras
%then $M \otimes^{\sigma {\rm h}} N$ is a unital operator space, from which it follows that
%its unital-subspace $M \otimes_{\rm w^*  h} N$ is also a unital operator space.

Since $\otimes_{\rm h}$ is projective, we have that $\otimes_{\rm w^*  h}$
is `weak*-injective'.  Thus if $X, Y$ are unital weak* closed
subspaces of $B(H)$ and $B(K)$ respectively, then
$X \otimes_{\rm w^*  h} Y$ is a unital weak* closed
subspace of $B(H) \otimes_{\rm w^*  h} B(K)$.
%D
However $B(H) \otimes_{\rm w^*  h} B(K)$ is a unital-subspace of
the unital dual operator space $B(H) \otimes^{\sigma {\rm h}} B(K)$.
The result is completed by an appeal to Theorem \ref{th2}, since
$X \otimes_{\rm w^*  h} Y$ is a dual operator space.
(We remark that one may in this case also construct an explicit  weak* continuous unital
 complete isometry from $B(H) \otimes_{\rm w^*  h} B(K)$
into a von Neumann algebra.)
\end{proof}

%\bigskip

%{\em Q
%DDWe close with two questions.  First, for which operator spaces $E$
%is $CB(E)$ a unital operator space?   This seems interesting in the
%light of recent results by
%Daws and Uygal that every completely contractive Banach
%algebra with predual, whose product is separately weak* continuos,
%is completely isometrically
%isomorphic to  a unital subalgebra of $CB(E)$ for a reflexive operator space $E$.
%  Second, if $X$
  %DD
 % is a dual operator system, then is every completely contractive weak*  continuous
% $u : X \to B(H)$ of the form $R v(\cdot)T$ for Hilbert space
%contractions $R, T$ and $v : X \to B(K)$ unital,
%weak*  continuous, and completely contractive?   In particular, the case $u : X \to \Cdb$ should be
%interesting to study.

%DD
Once one knows the spaces above are unital operator spaces, it is of
interest to compute their `noncommutative Shilov boundary', or
$C^*$-envelope \cite{Arv,BLM,P}.   It is easy to show for example
that the $C^*$-envelope of (even a module) Haagerup tensor product
of unital $C^*$-algebras is their full amalgamated free product
$C^*$-algebra (and we thank M. Junge for showing us a related fact).
Thus the $C^*$-envelope of $CB(M_n)$ is $M_n \,
*_{\Cdb} \, M_n$.

 \end{document}